\documentclass[11pt,twoside]{article}
\usepackage{graphicx}
\usepackage{amsmath}

\usepackage{amsmath}
\usepackage{amssymb}
\usepackage{amsfonts}
\usepackage[all]{xy}

\newtheorem{theorem}{Theorem}[section]
\newtheorem{proposition}[theorem]{Proposition}
\newtheorem{definition}[theorem]{Definition}
\newtheorem{conjecture}[theorem]{Conjecture}

\newcommand{\lto}{\longrightarrow}
\newcommand{\wt}{\widetilde}
\def\o#1{{\mathcal O}_{#1}}
\def\db#1{{{\mathcal D}^b}({#1})}
\def\H#1,#2,#3,#4{{\rm Hom}^{#1}_{#2}({#3}\:,\; {#4})}
\def\h#1,#2{{\rm Hom}({#1}\:,\; {#2})}

\title{\bf Derived Categories of Coherent Sheaves
}

\author{Alexei Bondal\vspace*{-0.5cm}\thanks{
Algebra Section, Steklov Mathematical Institute, Russian Academy
of Sciences,
 8 Gubkin str., GSP-1, Moscow 117966, Russia.
 E-mail: bondal@mi.ras.ru}\and Dmitri Orlov \vspace*{-0.5cm}\thanks{
Algebra Section, Steklov Mathematical Institute, Russian Academy
of Sciences,
 8 Gubkin str., GSP-1, Moscow 117966, Russia. E-mail: orlov@mi.ras.ru}}
\date{}

\begin{document}
\maketitle







\section{Introduction} \label{section 1}

This paper is devoted to studying the derived categories $\db {X}$ of coherent
sheaves on smooth algebraic varieties $X$ and on their noncommutative
counterparts. Derived categories of coherent sheaves proved to
contain the complete geometric information about varieties (in the sense of the
classical Italian school of algebraic geometry) as well as the
related homological algebra.

The situation when there exists a functor ${\db M}\lto {\db X}$
which is fully faithful is of special interest. We are convinced
that any example of such a functor is both algebraically and geometrically meaningful.

A particular case of a fully faithful functor is an equivalence of
derived categories ${\db M}\stackrel{\sim}{\lto}{\db X}$.

We show that smooth
projective varieties with ample canonical or anticanonical bundles
are uniquely determined by their derived categories. Hence the
derived equivalences between them boil down to autoequivalences.
We prove that
for such a variety the group of exact autoequivalences is the
semidirect product of the
group of automorphisms of the variety and the Picard group plus translations.

Equivalences and autoequivalences for the case of varieties with
non-ample (anti) canonical sheaf are now intensively studied.
The group of autoequivalences is believed to be closely related to
the holonomy group
of the mirror-symmetric family.

We give a criterion for a functor
between derived categories of coherent sheaves on two algebraic
varieties to be fully faithful. Roughly speaking, it is in
orthogonality of the images under the
functor of the structure sheaves of distinct closed points of the
variety. If a functor $\Phi: {\db M}\lto {\db X}$ is fully
faithful, then it induces  a so-called semiorthogonal
decomposition of ${\db X}$ into ${\db M}$ and its right orthogonal
category.

It turned out that derived categories have nice behavior under special
birational transformations like blow ups, flips and flops.
We describe a semiorthogonal decomposition of the derived category of the blow-up of a smooth variety $X$
in a smooth center $Y\subset X$. It contains one component isomorphic to $\db {X}$ and several components isomorphic to $\db {Y}$.

We also consider some flips and flops.
Examples support the conjecture that for any generalized flip $X\dashrightarrow X^+$
there exists a fully faithful functor $\db {X^+}\to \db {X}$ and it must be an equivalence for generalized flops.
This suggests the idea that the
minimal model program of the birational geometry can be viewed as
a `minimization' of the derived category $\db {X}$ in a
given birational class of $X$.

Then we widen the categorical approach to birational geometry by including in the picture some noncommutative varieties.
We propose to consider noncommutative desingularizations and formulate a conjecture generalizing the derived McKay correspondence.

We construct a semiorthogonal decomposition for the derived category of the complete
intersections of quadrics. It is related to classical questions of algebraic geometry, like 'quadratic complexes of lines',
and to noncommutative geometric version of Koszul quadratic duality.

The research described in this publication was made possible in part by Award No. RM1-2405-M0-02
of the U.S. Civilian Research and Development Foundation for the Independent States of the Former Soviet Union
(CRDF).

\section{Equivalences between derived categories}\label{section 2}

The first question that arises in studying algebraic varieties
from the point of view of derived categories is when varieties
have equivalent derived categories of coherent sheaves.
Examples of such equivalences for abelian varieties and K3
surfaces were constructed by Mukai \cite{Mu1}, \cite{Mu2},
Polishchuk \cite{Pol} and the second author in \cite{Or1}, \cite{Or2}.
See below on derived equivalences
for birational maps.

Yet we prove that a variety $X$ is uniquely determined by its
category ${\db{X}}$,
 if its anticanonical (Fano case) or canonical (general type case) sheaf
is ample.
To this end, we use only the graded (not triangulated) structure of the
category. By definition a {\sf graded category} is a pair $( {\cal
D}, T_{\cal D} )$ consisting of a category ${\cal D}$ (which we
always assume to be $k$-linear over a field $k$) and a fixed
equivalence $T_D : {\cal D}\lto {\cal D}$ , called translation
functor. For derived categories the translation
functor is defined to be the shift of grading in complexes.

Of crucial importance for exploring derived categories are
existence and properties of the Serre functor, defined in \cite{BK}.

\begin{definition}\label{SF}{\rm  \cite{BK} \cite{BO2}} Let ${\cal D}$ be a $k$-linear category
with finite--dimensional ${\rm Hom's}$. A covariant additive
functor $S: {\cal D}\to{\cal D}$ is called a {\sf Serre functor}
if it is an equivalence and there are given bi--functorial
isomorphisms for any $A,B\in{\cal D}$:
$$
\varphi_{A,B}: {\rm Hom}_{\cal D}(A\:,
\;B)\stackrel{\sim}{\lto}{\rm Hom}_{\cal D}(B\:, \;SA)^*.
$$

\end{definition}

A Serre functor in a category ${\cal D}$, if it exists, is unique
up to a graded natural isomorphism.

If $X$ is a smooth algebraic variety, $n={\rm dim}X$, then the
functor
$(\cdot)\otimes\omega_X[n]$
 is the Serre
functor in ${\db{X}}$. Thus, the Serre functor in ${\db{X}}$
can be viewed as a
categorical incarnation of the canonical sheaf $\omega_X$.

\begin{theorem}\label{rec}{\rm \cite{BO2}} Let $X$ be
a smooth irreducible projective variety with ample canonical or
anticanonical sheaf. If $\db{X}$ is equivalent
  as a graded category to $\db{X'}$ for some other smooth algebraic variety
$X'$, then $X$ is isomorphic to $X'$.
\end{theorem}

The idea of the proof is
that for varieties with ample canonical or anticanonical sheaf we
can recognize the skyscraper sheaves of closed points in $\db{X}$
by means of the Serre functor.
In this way we find the variety as a set. Then we reconstruct one
by one the set of line bundles, Zariski topology and the
structural sheaf of rings.

This theorem has a generalization to smooth orbifolds related to projective
varieties with mild singularities, as it was shown by Y. Kawamata \cite{Kaw1}.

Now consider the problem of computing the
group ${\rm Aut}{\db{X}}$ of exact (i.e. preserving triangulated structure)
autoequivalences of $\db{X}$ for an individual $X$.

\begin{theorem}\label{aut}{\rm \cite{BO2}}
Let $X$ be a smooth irreducible projective variety with
ample canonical or anticanonical sheaf. Then the group
of isomorphism classes of exact autoequivalences $\db{X} \to
 \db{X} $ is generated by the
automorphisms of the variety, twists by all invertible sheaves
and translations.
\end{theorem}

In the hypothesis of Theorem \ref{aut} the group ${\rm Aut}{\db{X}}$
is the semi-direct product of its subgroups
${\rm Pic} X \oplus {\bf Z}$ and ${\rm Aut} X$, ${\bf Z}$
being generated by the translation functor:
$$
{\rm Aut}{\db{X}}\cong {\rm Aut} X \ltimes ({\rm Pic} X \oplus {\bf Z}).
$$

\section{Fully faithful functors and semiorthogonal
decompositions}\label{section 3}

An equivalence is a particular instance of a fully faithful functor.
This is a functor $F:{\cal C}\to {\cal
D}$ which for any pair of objects $X,Y\in {\cal C}$ induces an
isomorphism ${\h X, Y}\simeq {\h FX, {FY}}$. This notion is especially
useful in the context of triangulated categories.

If a functor $\Phi: {\db M}\lto {\db X}$ is fully faithful, then
it induces  a so-called semiorthogonal decomposition of ${\db X}$
into ${\db M}$ and its right orthogonal.


Let ${\cal B}$ be a full triangulated subcategory of a
triangulated category ${\cal D}$. The {\sf right orthogonal} to
${\cal B}$ is the full subcategory ${\cal B}^{\perp}\subset {\cal
D}$ consisting of the objects $C$ such that ${\h B, C}=0$ for all
$B\in{\cal B}$. The {\sf left orthogonal} ${}^{\perp}{\cal B}$ is
defined analogously.
The categories ${}^{\perp}{\cal B}$ and ${\cal B}^{\perp}$ are
also triangulated.






\begin{definition}\label{sd} {\rm \cite{BK}}
A sequence of triangulated subcategories $({\cal B}_0,..., {\cal
B}_n)$ in a triangulated category ${\cal D}$ is said to be {\sf
semiorthogonal} if ${\cal B}_j\subset {\cal
B}^{\perp}_i$ whenever $0\le j<i\le n$. If a semiorthogonal sequence
generates ${\cal D}$ as a triangulated category, then
we call it by {\sf
semiorthogonal decomposition} of the category ${\cal D}$ and
denote this as follows:
$$
{\cal D}=\Bigl\langle{\cal B}_0,....,{\cal B}_n\Bigl\rangle.
$$
\end{definition}

Examples of semiorthogonal decompositions are provided by exceptional
sequences of objects \cite{B}. These arise when all
${\cal B}_i$'s are equivalent to the derived categories of finite
dimensional  vector spaces ${\db {k-{\rm mod}}}$. Objects which
correspond to the 1-dimensional vector space under a fully
faithful functor $F:{\db {k-{\rm mod}}}\to {\mathcal D}$ can be intrinsically
defined as  {\em exceptional}, i.e. those satisfying the
conditions ${\rm Hom}^i(E\:,\;E)=0$, when $i\ne0$, and ${\rm
Hom}^0(E\:,\;E)=k$. There is a natural action of the braid group on exceptional sequences \cite{B}
and, under some conditions, on semiorthogonal sequences of subcategories in a triangulated category \cite{BK}.

We propose to consider the derived category of coherent sheaves as an
analogue of the motive of a variety, and semiorthogonal
decompositions as a tool for  simplification of this `motive' similar
to splitting by projectors in Grothendieck motivic theory.


Let $X$ and $M$ be smooth algebraic varieties of dimension $n$ and $m$
respectively and $E$ an object in ${\db {X \times M}}$.
Denote by $p$ and $\pi$ the projections of $M \times X$
to $M$ and $X$ respectively.
With $E$ one can associate the functor $ \Phi_{E} : {\db M}\lto {\db X}$
 defined by the formula:
$$
\Phi_{E}(\cdot):={\mathbf R}\pi_* (E \stackrel{\mathbf L}{\otimes}
p^*(\cdot)).
$$

It happens that any fully faithful functor is of this form.

\begin{theorem}\label{main} {\rm \cite{Or1}} Let $F:\db{M}\to \db{X}$ be an exact fully faithful functor,
where $M$ and $X$ are smooth projective
varieties.
Then there exists a unique up to isomorphism object $E\in\db{M\times X}$ such
 that $F$ is isomorphic to the functor $\Phi_E$.
\end{theorem}

The assumption on existence of the right adjoint to $F$, which was originally in \cite{Or1}, can be removed in view
of saturatedness of $\db {M}$ due to \cite{BK}, \cite{BVdB}.

This theorem is in conjunction with the following criterion.

\begin{theorem}\label{mai}{\rm \cite{BO1}}
Let $M$ and $X$ be smooth algebraic varieties and
 $E\in{\db {M\times X}}$. Then $\Phi_{E}$ is  fully faithful functor
if and only if the following orthogonality conditions are verified:
$$
\begin{array}{lll}
i) & {\H i, X, \Phi_E({\cal O}_{t_1}), {\Phi_{E}({\cal O}_{t_2})}} = 0 & \qquad \mbox{for every }\: i\;\mbox{ and } t_1\ne t_2.\\
&&\\
ii) & {\H 0, X, \Phi_E({\cal O}_t), {\Phi_E({\cal O}_t)}} = k,&\\
&&\\
& {\H i, X, \Phi_E({\cal O}_t), {\Phi_E({\cal O}_t)}} = 0 , & \qquad \mbox{ for }i\notin [0, dim M].
\end{array}
$$
Here $t$, $t_{1}$, $t_{2}$ stand for closed points in $M$ and
${\cal O}_{t_{i}}$ for the skyscraper sheaves.
\end{theorem}

The criterion is a particular manifestation of the following important
principle:
{\sf
suppose $M$ is realized as an appropriate moduli space of
pairwise homologically orthogonal objects in a triangulated
category ${\cal D}$ taken `from the real life',  then one can expect a sheaf
of finite (noncommutative) algebras ${\cal A}_M$ over ${\cal O}_M$
and a fully faithful functor from the derived category
$\db {coh({\mathcal A}_M)}$ of coherent modules over ${\cal A}_M$ to ${\cal D}$.
}

There are also strong indications that this principle should have
a generalization, at the price of considering noncommutative DG moduli
spaces, to the case when the orthogonality condition is
dropped.

\section{Derived categories and birational geometry}
\label{section 5}

Behavior of derived categories under birational
transformations shows that they can serve as a useful tool in comprehending
various phenomena of birational geometry and play possibly the key role in
realizing the minimal model program.

First, we give a description of the derived category of
the blow-up of a variety with smooth center in terms of the
categories of the variety and of the center.
Let $Y$ be a smooth subvariety of codimension $r$ in a smooth
algebraic variety $X$. Denote $\widetilde X$ the smooth algebraic
variety obtained by blowing up $X$ along the center $Y$.
There exists a fibred square:
$$
\begin{array}{ccc}
\widetilde Y&\stackrel{j }{\longrightarrow}& \widetilde X\\
\llap{\footnotesize $p$} \downarrow &&\llap{$\pi$} \downarrow\\
Y&\stackrel{i}{\longrightarrow}& X
\end{array}
$$
where $i$ and $j$ are smooth embeddings, and $p:
\widetilde Y\to Y$ is the projective fibration of the exceptional
divisor $\widetilde Y$ in $\widetilde X$ over the center $Y$.
Recall that $\widetilde Y\cong {\mathbb P}(N_{X/Y})$ is the
projective normal bundle.
Denote by $\o{\widetilde Y}(1)$ the relative Grothendieck sheaf.

\begin{proposition}\label{blow} {\rm  \cite{or}}
Let  $L$ be any invertible sheaf on $\widetilde Y$.
The  functors
$$
{\mathbf L}\pi^*: \db {X}\lto \db {\widetilde X},
$$
$$
{\mathbf R}j_*(L\otimes p^*(\cdot)): \db {Y}\lto \db {\widetilde X}
$$
are fully faithful.
\end{proposition}

Denote by $D(X)$ the full subcategory of $\db{\widetilde X}$ which is the image of
$\db{X}$ with respect to the functor ${\mathbf L}\pi^*$ and by $D(Y)_k$ the full subcategories of
$\db{\widetilde X}$ which are the images of $\db{Y}$ with respect to the functors ${\mathbf R}j_*(\o{\widetilde Y}(k)\otimes p^*(\cdot))$.

\begin{theorem}\label{full}{\rm \cite{or}}{\rm \cite{BO1}}
We have the semiorthogonal decomposition of the category of the blow-up:
$$
\db{\widetilde X}=\Bigl< D(Y)_{-r+1}, ... , D(Y)_{-1}, D(X) \Bigl>.
$$

\end{theorem}

Now we consider the behavior of the derived categories of coherent sheaves with respect to the special birational transformations
called flips and flops.

Let $Y$ be a smooth subvariety of a smooth algebraic variety $X$
such that $Y\cong {\mathbb P}^{k}$ and $N_{X/Y}\cong
{\mathcal O}(-1)^{\oplus(l+1)}$ with $l\le k$.

If $\widetilde X$ is the blow-up of $X$ along $Y$, then  the exceptional
divisor $\widetilde Y\cong {\mathbb P}^k\times {\mathbb P}^l$ is  the
product of projective spaces. We can blow down $\widetilde X$ in
such a way that $\widetilde Y$ projects to the second component ${\mathbb P}^l$
of the product. As a result we obtain a smooth  variety $X^+$, which for
simplicity we assume to be algebraic, with subvariety
$Y^+\cong {\mathbb P}^l$. This is depicted in the following diagram:
$$
\xymatrix{
&\wt{Y}\ar[d]^j \ar[dl]_p\ar[dr]^{p^+}&\\
Y\ar[d]_i& \wt{X}\ar[dl]_(.4){\pi}
\ar[dr]^(.4){\pi^+}&
Y^{+}\ar@<-0.8ex>[d]^{i^+}\\
X\ar@{-->}[rr]^{fl}&&X^{+}}
$$
 The birational map $X\dashrightarrow X^+$ is the simplest instance of flip, for $l\le k$.
If $l=k$, this is a flop.
\begin{theorem}\label{fl}{\rm \cite{BO1}}
In the above notations, the functor ${\mathbf R}\pi_* {\mathbf L}\pi^{+*}: \db {X^+}\lto \db {X}$
is fully faithful for $l\le k$. If $l=k$, it is an equivalence.
\end{theorem}

This theorem has an obvious generalization to the case when $Y$ is isomorphic to
the projectivization of a vector bundle $E$ of rank $k$ on a smooth variety $W$, $q: Y\lto W$,  and
$N_{X/Y} =q^* F \otimes {\mathcal O}_{E}(-1)$ where $F$ is a vector bundle on $W$ of rank $l\le k$.
Then the blow-up with a smooth center can be viewed as the particular case of this flip when $Y$ is a divisor in $X$.
 Kawamata \cite{Kaw1} generalized the theorem to those flips between smooth
 orbifolds which are elementary (Morse type) cobordisms in the theory of birational cobordisms due to Wlodarczyk {\em et al.} \cite{Wl}, \cite{AKMW}.

Let $X$ and $X^+$ be  smooth projective varieties.
A birational map $ X\stackrel{fl}{\dashrightarrow} X^+$ will be called a {\sf generalized flip} if for some (and consequently for any) its
smooth resolution
$$
\xymatrix{
& \wt{X}\ar[dl]_(.4){\pi}
\ar[dr]^(.4){\pi^+}&\\
X\ar@{-->}[rr]^{fl}&&X^{+}}
$$
the difference $D=\pi^* K_X - \pi^{+*} K_{X^+}$ between the pull-backs of the canonical
divisors is an effective divisor on $\wt{X}$. The particular case when $D=0$ is called {\sf generalized flop}.

Theorem \ref{fl} together with calculations of 3-dimensional flops with
centers in  $(-2)$-curves \cite{BO1}
lead us to the following conjecture.

\begin{conjecture} For any generalized flip $X\stackrel{fl}{\dashrightarrow} X^+$
there is an exact fully faithful functor
$F: {\db{X^+}}\lto {\db{X}}$. It is  an equivalence for generalized flops.
\end{conjecture}


This  conjecture was
recently proved in dimension 3 by T. Bridgeland \cite{Br}.

The functor ${\mathbf R}\pi_* {\mathbf L}\pi^{+*}: \db {X^+}\lto \db {X}$ is not always fully faithful
under conditions of the conjecture, but we expect that it is such when $\wt{X}$ is replaced by the fibred product of $X$ and $X^+$ over
some common singular contraction of $X$ and $X^+$. Namikawa proved
 that this is the case for Mukai symplectic flops \cite{Na}.

A fully faithful functor $\db {X^+}\lto \db {X}$ induces a semiorthogonal decomposition of $\db {X}$ into $\db {X^+}$ and its right
orthogonal (which is trivial for flops).
Hence, passing from $X$ to $X^+$ for generalized flips has the categorical meaning
of breaking off semiorthogonal summands from the derived category. This suggests the idea that the minimal model program
 of birational geometry should be
interpreted as a minimization of the derived category $\db {X}$ in a given birational class of $X$.
Promisingly, chances are that the very existence of flips can be achieved by constructing $X^+$
as an appropriate moduli space of objects in $\db {X}$,
in accordance with the principle of the previous section (this is done by T. Bridgeland for flops in dimension 3 \cite{Br}).

\section{Noncommutative resolutions of singularities}
\label{section 55}

In this section we will give a  perspective
for categorical interpretation of the minimal model program by
enriching  the landscape with
the derived categories of noncommutative varieties.

Let $\pi :{\wt X}\to X$ be a proper birational morphism, where $X$ has
rational singularities. Then ${\mathbf R}\pi_*:\db {\wt X}\to \db {X}$ identifies
$\db {X}$ with the quotient of $\db {\wt X}$ by the kernel of ${\mathbf R}\pi_*$.
For this reason, let us call by a {\em categorical desingularization} of
a triangulated category ${\mathcal D}$ a pair $({\mathcal C}, {\mathcal K})$
consisting of an abelian category ${\mathcal C}$ of finite homological dimension and of
${\mathcal K}$, a thick subcategory in $\db {\mathcal C}$
such that ${\mathcal D}=\db {\mathcal C}/{\mathcal K}$.
We expect that for ${\mathcal D}=\db {X}$ there exists a minimal desingularization,
i.e. such one that $\db {\mathcal C}$ has a fully faithful embedding
in  $\db {{\mathcal C}'}$ for any other categorical desingularization $({\mathcal C}', {\mathcal K}')$ of ${\mathcal D}$.
Such a desingularization is unique up to derived equivalence of ${\mathcal C}$.

For the derived categories of singular varieties one may hope to find the minimal
desingularizations in the spirit of noncommutative geometry.

Let $X$ be an algebraic variety. We call by
{\sf noncommutative (birational) desingularization} of $X$ a pair $(p, {\mathcal A})$
consisting of a proper birational morphism $p:Y\to X$ and an algebra
 ${\mathcal A}={\mathcal E}nd({\mathcal F})$ on $Y$, the sheaf of local
endomorphisms of a torsion free coherent
${\mathcal O}_Y$-module ${\mathcal F}$, such that the abelian category of coherent
${\mathcal A}$-modules has finite homological dimension.

When $f: Y\to X$ is a morphism from a smooth $Y$ onto an affine $X$ with fibres of dimension $\le 1$ and ${\mathbf R}f_*({\mathcal O}_Y)={\mathcal O}_X$,
 M. Van den Bergh \cite{VdB} has recently constructed a noncommutative desingularization of $X$, which is derived equivalent to $\db {Y}$.


\begin{conjecture}
Let $X$ has canonical singularities and $q:Y\to X$ a finite morphism with smooth $Y$. Then
the pair $(id_{X}, {\mathcal E}nd(q_*{\mathcal O}_Y))$ is a minimal desingularization of $X$.
\end{conjecture}
In particular, we expect that $\db {coh({\mathcal E}nd(q_*{\mathcal O}_Y))}$ has a fully faithful functor into $\db {\wt X}$ for any smooth (commutative)
resolution of $X$. Moreover, if the resolution is crepant then the functor has to be an equivalence.

Let $X$ be the quotient of a smooth $Y$ by an action of a finite group $G$. If the locus of the points in $Y$ with nontrivial stabilizer in $G$
has codimension $\ge 2$, then the category of coherent ${\mathcal E}nd(q_*{\mathcal O}_Y)$-modules is equivalent to the category
of $G$-equivariant coherent sheaves on $Y$.
Therefore the conjecture is a generalization of
the derived McKay correspondence due to Bridgeland-King-Reid \cite{BKR}.

\section{Complete intersection of quadrics and noncommutative geometry}
\label{section 6}

This section is related to the previous one by Grothendieck slogan that projective geometry is a part
of theory of singularities.

Let $X$ be a smooth intersection of two projective
quadrics of even dimension $d$ over an algebraically closed field of characteristic zero.
 It appears that if we consider the hyperelliptic curve $C$ which
is the double cover of ${\mathbb P}^1$ that parameterizes the
pencil of quadrics, with ramification in the points corresponding
to degenerate quadrics, then ${\db C}$ is embedded in ${\db X}$ as
a full subcategory \cite{BO1}.
This gives a categorical explanation for the classical description
of moduli spaces of semistable bundles on the curve $C$ as moduli spaces of
(complexes of) coherent sheaves on $X$.

The orthogonal to ${\db C}$ in ${\db X}$ is
decomposed into an exceptional sequence (of line bundles ).
More precisely, we have a semiorthogonal decomposition
\begin{equation}\label{2dim}
\db{X}\ =\ \Bigl\langle{\cal O}_X(-d+3),\ldots,{\cal O}_X,\db{C}\Bigl\rangle .
\end{equation}

When a greater number of
quadrics is intersected, objects of noncommutative geometry naturally show up:
instead of coherent sheaves on hyperelliptic curves we must consider modules
over a sheaf of
 noncommutative algebras. More
about noncommutative geometry is in the talk of T. Stafford at
this Congress.

Consider a system of $m$ quadrics in ${\mathbb P}(V)$,
i.e. a linear embedding $U\stackrel{\phi}{\hookrightarrow} S^2 V^*$, where $\dim U=m$, $\dim V=n$, $2m\le n$.
Let $X$, the complete intersection of the quadrics, be a smooth subvariety in ${\mathbb P}(V)$ of dimension $n-m-1$.
Let $A=\mathop{\oplus}\limits_{i\ge 0} H^0(X, {\mathcal O}(i))$ be the coordinate ring of $X$. This graded
quadratic algebra is Koszul due to Tate \cite{Ta}. The quadratic dual algebra $B=A^!$ is the generalized homogeneous
Clifford algebra. It is generated in degree 1 by the space $V$, the relations being given
by the kernel of the dual to $\phi$ map $S^2V\to U^*$, viewed as a subspace in $V\otimes V$.
The center of $B$ is generated by $U^*$ (a subspace of quadratic elements in $B$) and an element $d$,
which satisfies the equation $d^2=f$ where $f$ is the equation of the locus of degenerate quadrics in $U$.
Algebra $B$ is finite over the central subalgebra $S=S^{\bullet}U^*$. The Veronese subalgebra $B_{ev}=\oplus B_{2i}$
is finite over the Veronese subalgebra $S_{ev}=\oplus S^{2i}U^*$.
Since ${\bf Proj}\ S_{ev}$ is isomorphic to $\mathbb{P}(U)$, the
sheafification
of $B_{ev}$ over  ${\bf Proj}\ S_{ev}$ is a sheaf ${\mathcal B}$ of finite
algebras over ${\mathcal O}_{\mathbb{P}(U)}$.
Consider the derived category $\db {coh({\mathcal B})}$ of coherent right ${\mathcal B}$-modules.

\begin{theorem}\label{pq}
Let $X$ be the smooth intersection of $m$ quadrics in ${\mathbb P}^{n-1}$, $2m\le n$.
Then there exists a fully faithful functor $\db {coh({\mathcal B})}\to \db {X}$. Moreover,
\begin{itemize}
\item[{\rm (i)}]
if $2m<n$, we have a semiorthogonal decomposition
$$
\db{X}\ =\ \Bigl\langle{\cal O}_X(-n+2m+1),\ldots,{\cal O}_X,\db{coh({\mathcal B})}\Bigl\rangle,
$$
\item[{\rm (ii)}] if $2m=n$, there is an equivalence
$\db {coh({\mathcal B})}\stackrel{\sim}{\to} \db {X}.$
\end{itemize}
\end{theorem}

For $m=0$, i.e. when there is no quadrics, the theorem coincides with
Beilinson's description of the derived category of the
projective space \cite{Be}. For $m=1$, this is Kapranov's description
of the derived category of the quadric \cite{Kap}.

For odd $n$, the element $d$ generates the center of ${\mathcal B}$ over ${\mathcal O}_{{\mathbb P}(U)}$.
Hence the spectrum of the center of ${\mathcal B}$ is a ramified double cover $Y$ over ${\mathbb P}(U)$.
Also ${\mathcal B}$ yields a coherent sheaf of algebras ${\mathcal B}'$ over $Y$,
such that $coh({\mathcal B}')$ is equivalent to $coh({\mathcal B})$. For the
above case of two even dimensional quadrics,
${\mathcal B}'$ is an Azumaya algebra over $Y=C$. Since Brauer group of $Y$ (taken over an
algebraically closed field of characteristic zero) is trivial, the category $coh({\mathcal B}')$ is equivalent to $coh({\mathcal O}_Y)$.
Hence (\ref{2dim}) follows from the theorem.

Furthermore, when $X$ is a K3 surface, the smooth intersection of 3 quadrics in ${\mathbb P}^5$, then the double cover $Y$ is also
a K3 surface, but ${\mathcal B}'$ is in general a nontrivial Azumaya algebra over $Y$. The theorem states an equivalence
$\db {X}\simeq \db{coh({\mathcal B}')}$.

This theorem illustrates the principle from section 3. The fully faithful functor is related to the moduli space
of vector bundles on $X$, which are the restrictions to $X$ of the spinor bundles on the quadrics.
The (commutative) moduli space involved is either ${\mathbb P}(U)$ or $Y$, depending on parity of $n$.

Algebraically, the fully faithful
functor in the theorem is given by an appropriate version of Koszul duality.
Theorem \ref{pq} has a generalization to a class of Koszul Gorenstein
algebras, which includes the coordinate rings of superprojective spaces.

\end{document}